\documentclass[12pt,twoside]{article} 
\usepackage[english]{babel}
\usepackage[latin1]{inputenc} 
\usepackage{amsmath}
\usepackage{amssymb,amsfonts}
\usepackage{graphicx}                   
\usepackage{times,amssymb,amscd}

\newcommand{\bE}{\mathbf{E}}
\newcommand{\bG}{\mathbf{G}}
\newcommand{\bH}{\mathbf{H}}

\newcommand{\bL}{\mathbf{L}}

\newcommand{\bR}{\mathbf{R}}
\newcommand{\bS}{\mathbf{S}}
\newcommand{\bV}{\mathbf{V}}

\newcommand{\bT}{\mathbf{T}}

\newcommand{\BV}{\boldsymbol{V}}

\newcommand{\cP}{\mathcal{P}}

\newcommand{\cT}{\mathcal{T}}

\newcommand{\cH}{\mathcal{H}}

\newcommand{\EUC}{\bE^3}

\newcommand{\SXR}{\bS^2\!\times\!\bR}
\newcommand{\HXR}{\bH^2\!\times\!\bR}
\newcommand{\SLR}{\widetilde{\bS\bL_2\bR}}
\newcommand{\NIL}{\mathbf{Nil}}
\newcommand{\SOL}{\mathbf{Sol}}

\newtheorem{Definition}{Definition}[section]
\newtheorem{Remark}{Remark}[section]

\begin{document}
\pagestyle{myheadings}
\markboth{\centerline{Jen\H o Szirmai}}
{Regular prism tilings in $\SLR$ space}
\title
{Regular prism tilings in $\SLR$ space }

\author{Jen\H o Szirmai \\
\normalsize Budapest University of Technology and \\ 
\normalsize Economics Institute of Mathematics, \\
\normalsize Department of Geometry \\
\normalsize Budapest, P. O. Box: 91, H-1521 \\
\normalsize szirmai@math.bme.hu
\date{\normalsize{\today}}}


\maketitle
\begin{abstract}

$\SLR$ geometry is one of the eight 3-dimensional  
Thurston geometries, it can be derived from the 3-dimensional Lie group of all $2\times 2$ real matrices with determinant one.

Our aim is to describe and visualize the {\it regular infinite (torus-like) or bounded} $p$-gonal prism tilings in $\SLR$ space. 
For this purpose we introduce the notion of the infinite and bounded prisms, 
prove that there exist infinite many regular infinite $p$-gonal face-to-face prism tilings $\cT^i_p(q)$ and 
infinitely many regular (bounded) $p$-gonal non-face-to-face $\SLR$ prism tilings $\cT_p(q)$ for parameters $p \ge 3$ 
where $ \frac{2p}{p-2} < q \in \mathbb{N}$.  
Moreover, we develope a method to determine the data of the space filling regular infinite 
and bounded prism tilings. We apply the above procedure to $\cT^i_3(q)$ and $\cT_3(q)$ where $6< q \in \mathbb{N}$ and visualize them 
and the corresponding tilings.

E. {Moln\'ar} showed, that the homogeneous 3-spaces
have a unified interpretation in the projective 3-space $\mathcal{P}^3(\bV^4,\BV_4, \mathbf{R})$. 
In our work we will use this projective model of $\SLR$
geometry and in this manner the  prisms and prism tilings can be visualized on the Euclidean screen of computer.

\end{abstract}

\newtheorem{theorem}{Theorem}[section]
\newtheorem{corollary}[theorem]{Corollary}
\newtheorem{conjecture}[theorem]{Conjecture}
\newtheorem{lemma}[theorem]{Lemma}
\newtheorem{exmple}[theorem]{Example}
\newenvironment{definition}{\begin{defn}\normalfont}{\end{defn}}
\newenvironment{remark}{\begin{rmrk}\normalfont}{\end{rmrk}}
\newenvironment{example}{\begin{exmple}\normalfont}{\end{exmple}}
\newenvironment{acknowledgement}{Acknowledgement}



\section{On $\SLR$ geometry}

The real $ 2\times 2$ matrices $\begin{pmatrix}
         d&b \\
         c&a \\
         \end{pmatrix}$ with unit determinant $ad-bc=1$ 
constitute a Lie transformation group by the usual product operation, taken to act on row matrices as on point coordinates on the right as follows
\begin{equation}
(z^0,z^1)\begin{pmatrix}
         d&b \\
         c&a \\
         \end{pmatrix}=(z^0d+z^1c, z^0 b+z^1a)=(w^0,w^1). \tag{1.1}
\end{equation}
This group is a $3$-dimensional manifold, because of its $3$ independent real coordinates and with its usual neighbourhood topology (\cite{S}, \cite{T}).
In order to model the above structure on the projective space $\cP^3$ (see \cite{M97}) we introduce the new projective coordinates $(x^0,x^1,x^2,x^3)$ where
\begin{equation}
a:=x^0+x^3, \ b:=x^1+x^2, \ c:=-x^1+x^2, \ d:=x^0-x^3, \notag
\end{equation}
with positive equivalence as a projective freedom. 
Then it follows, that
\begin{equation}
0>bc-ad=-x^0x^0-x^1x^1+x^2x^2+x^3x^3 \tag{1.2}
\end{equation}  
describes the interior of the above one-sheeted hyperboloid solid $\cH$ in the usual Euclidean coordinate simplex with the origin
$E_0(1;0;0;0)$ and the ideal points of the axes $E_1^\infty(0;1;0;0)$, $E_2^\infty(0;0;1;0)$, $E_3^\infty(0;0;0;1)$.
We consider the collineation group ${\bf G}_*$ which acts on the projective space $\cP^3$  and preserves a polarity i.e. a scalar product of signature
$(- - + +)$, this group leave the one sheeted hyperboloid solid $\cH$ invariant. 
We have to choice a appropriate subgroup $\mathbf{G}$ of $\mathbf{G}_*$ as isometry group, then the universal covering space 
$\widetilde{\cH}$ of $\cH$ will be the hyperboloid model of $\SLR$ (see \cite{M97}).

The specific isometry $\bS$ is an one parameter group given by the matrices $(s_i^j(\phi))$:
\begin{equation}
\begin{gathered} \bS(\phi):~(s_i^j(\phi))=
\begin{pmatrix}
\cos{\phi}&\sin{\phi}&0&0 \\
-\sin{\phi}&\cos{\phi}&0&0 \\
0&0&\cos{\phi}&-\sin{\phi} \\
0&0&\sin{\phi}&\cos{\phi}
\end{pmatrix}
\end{gathered} \tag{1.3}
\end{equation}
The elements of $\bS$ are the so-called {\it fibre translations}. We obtain an unique fibre line to each $X(x^0;x^1;x^2;x^3) \in \widetilde{\cH}$ 
as the orbit by right action of $\bS$ on $X$. The coordinates of points lying on the fibre line through $X$ can be expressed 
as the images of $X$ by $\bS(\phi)$:
\begin{equation}
\begin{gathered}
(x^0;x^1;x^2;x^3) \stackrel{\bS(\phi)}{\longrightarrow} {(x^0 \cos{\phi}-x^1 \sin{\phi}; x^0 \sin{\phi} + x^1 \cos{\phi};} \\ {x^2 \cos{\phi} + x^3 \sin{\phi};-x^2 \sin{\phi}+
x^3 \cos{\phi})}.
\end{gathered} \tag{1.4}
\end{equation}
The points of a fibre line throught 
$X$ by usual inhomogeneous Euclidean coordinates $x=\frac{x^1}{x^0}$, $y=\frac{x^2}{x^0}$, $z=\frac{x^3}{x^0}$, $x^0\ne 0$ are given by
\begin{equation}
\begin{gathered}
(1;x;y;z) \stackrel{\bS(\phi)}{\longrightarrow} {\Big( 1; \frac{x+\tan{\phi}}{1-x \tan{\phi}}; \frac{y+z \tan{\phi}}{1-x \tan{\phi}};
\frac{z - y \tan{\phi}}{1-x \tan{\phi}}\Big)}.
\end{gathered} \tag{1.5}
\end{equation}
\begin{figure}[ht]
\centering
\includegraphics[width=10cm]{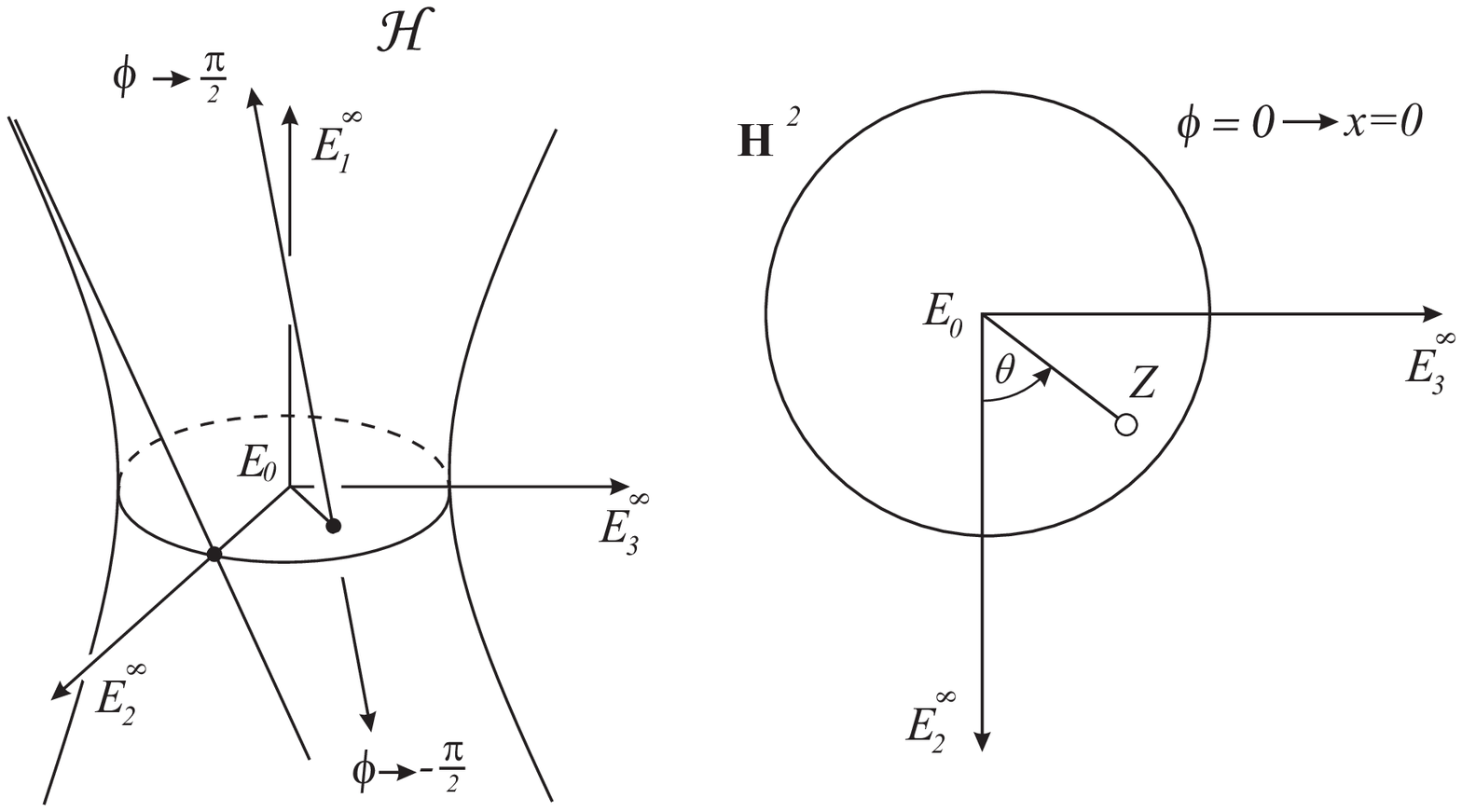}
\caption{}
\label{}
\end{figure}
The $\pi$ periodicity of the above maps can be seen from the formulas (1.4) and (1.5) e.g. if $-\frac{\pi}{2} < \phi < \frac{\pi}{2}$ then 
$-\infty < x < \infty$.
The elements of the isometry group of $\SLR$ in the above basis can be described by the matrix $(a_i^j)$ (see \cite{MSz})
\begin{equation}
\begin{gathered} (a_i^j)=
\begin{pmatrix}
a_0^0&a_0^1&a_0^2&a_0^3 \\
\mp a_0^1&\pm a_0^0&\pm a_0^3&\mp a_0^2 \\
a_2^0&a_2^1&a_2^2&a_2^3 \\
\pm a_2^1&\mp a_2^0&\mp a_2^3&\pm a_2^2 \\
\end{pmatrix} \ \ \text{where} \\
-(a_0^0)^2-(a_0^1)^2+(a_0^2)^2+(a_0^3)^2=-1, \ \ -(a_2^0)^2-(a_2^1)^2+(a_2^2)^2+(a_2^3)^2=1, \\
-a_0^0a_2^0-a_0^1a_2^1+a_0^2a_2^2+a_0^3a_2^3=0=-a_0^0a_2^1+a_0^1a_2^0-a_0^2a_2^3+a_0^3a_2^2.
\end{gathered} \tag{1.6}
\end{equation}
We define the {\it translation group} $\bG_T$ as a subgroup of $\SLR$ isometry group acting transitively on the points of $\widetilde{\cH}$ and 
mapping the origin $E_0(1;0;0;0)$ onto $X(x^0;x^1;x^2;x^3;)$. These isometries and their inverses (up to a positive determinant factor) can be given by the following $(t_i^j)$ and
$T_j^k=(t_i^j)^{-1}$ matrices:
\begin{equation}
\begin{gathered} \bT:~(t_i^j)=
\begin{pmatrix}
x^0&x^1&x^2&x^3 \\
-x^1&x^0&x^3&-x^2 \\
x^2&x^3&x^0&x^1 \\
x^3&-x^2&-x^1&x^0
\end{pmatrix},\\ \bT^{-1}:~(T_j^k)=
\begin{pmatrix}
x^0&-x^1&-x^2&-x^3 \\
x^1&x^0&-x^3&x^2 \\
-x^2&-x^3&x^0&-x^1 \\
-x^3&x^2&x^1&x^0
\end{pmatrix}.
\end{gathered} \tag{1.7}
\end{equation}
The rotation about the fibre line through the origin $E_0(1;0;0;0)$ by angle $\omega$ $(-\pi<\omega\le \pi)$ can be expressed by the following matrix
(see (1.8) and \cite{M97})
\begin{equation}
\begin{gathered} \bR_{E_O}(\omega):~(r_i^j(E_0,\omega))=
\begin{pmatrix}
1&0&0&0 \\
0&1&0&0 \\
0&0&\cos{\omega}&\sin{\omega} \\
0&0&-\sin{\omega}&\cos{\omega}
\end{pmatrix},
\end{gathered} \tag{1.8}
\end{equation}
and the rotation $\bR_X(\omega)$ about the fibre line through $X(x^0;x^1;x^2;x^3)$ by angle $\omega$ can be derived by formulas (1.7) and (1.8): 
\begin{equation}
\begin{gathered} \bR_X(\omega)=\bT^{-1} \bR_{E_O} (\omega) \bT:~(r_i^j(X,\omega))=\\
\begin{pmatrix}
x^0&-x^1&-x^2&-x^3 \\
x^1&x^0&-x^3&x^2 \\
-x^2&-x^3&x^0&-x^1 \\
-x^3&x^2&x^1&x^0
\end{pmatrix}
\begin{pmatrix}
1&0&0&0 \\
0&1&0&0 \\
0&0&\cos{\omega}&\sin{\omega} \\
0&0&-\sin{\omega}&\cos{\omega}
\end{pmatrix}
\begin{pmatrix}
x^0&x^1&x^2&x^3 \\
-x^1&x^0&x^3&-x^2 \\
x^2&x^3&x^0&x^1 \\
x^3&-x^2&-x^1&x^0
\end{pmatrix}.
\end{gathered} \tag{1.9}
\end{equation}
Horizontal intersection of the hyperboloid solid $\cH$ e.g. with the plane $E_0^\infty E_2^\infty E_3^\infty$ provide the Beltrami-Cayley-Klein 
model of the hyperbolic plane $\bH^2$ that is called {\it base plane} of the model $\widetilde{\cH}=\SLR$. 
The fibre through $X$ intersects the $z^1=x=0$ base plane in a trace point
\begin{equation}
\begin{gathered}
Z(z^0=x^0 x^0+x^1x^1; z^1=0; z^2=x^0x^2-x^1x^3;z^3=x^0x^3+x^1x^2).
\end{gathered} \tag{1.10}
\end{equation}
We introduce a so-called hyperboloid parametrization by \cite{M97} as follows
\begin{equation}
\begin{gathered}
x^0=\cosh{r} \cos{\phi}, \\
x^1=\cosh{r} \sin{\phi}, \\
x^2=\sinh{r} \cos{(\theta-\phi)}, \\
x^3=\sinh{r} \sin{(\theta-\phi)}, 
\end{gathered} \tag{1.11}
\end{equation}
where $(r,\theta)$ are the polar coordinates of the base plane and $\phi$ is just the fibre coordinate. We note that
$$-x^0x^0-x^1x^1+x^2x^2+x^3x^3=-\cosh^2{r}+\sinh^2{r}=-1<0.$$
The inhomogeneous coordinates corresponding to (1.11), that play an important role in later visualization of the prism tilings in $\EUC$,
are given by
\begin{equation}
\begin{gathered}
x=\frac{x^1}{x^0}=\tan{\phi}, \\
y=\frac{x^2}{x^0}=\tanh{r} \frac{\cos{(\theta-\phi)}}{\cos{\phi}}, \\
z=\frac{x^3}{x^0}=\tanh{r} \frac{\sin{(\theta-\phi)}}{\cos{\phi}}.
\end{gathered} \tag{1.12}
\end{equation}
\section{Prisms and prism tilings in $\SLR$ space}
After having investigated the prisms and prism-like tilings in $\SXR$ and $\HXR$ spaces
(see \cite{Sz10-2} and \cite{Sz10-3}) we consider the analogous problem in $\SLR$ space from among
the eight Thurston geometries. 
\begin{Definition}
Let $\cP^i$ be a $\SLR$ infinite solid that is bounded by one-sheeted hyperboloid surfaces of the model space generated by neighbouring ,,side fibre lines" passing through the 
vertices of a $p$-gon ($\cP^b$) lying in the ,,hyperbolic base plane".
The images of solids $\cP^i$ by $\SLR$ isometry are called infinite (or torus-like)  $p$-sided $\SLR$ prisms. 
\end{Definition}
The cammon part of $\cP^i$ with the hyperbolic base plane is the {\it base figure} of $\cP^i$ that is denoted by $\cP$ and its vertices coincide
with the vertices of $\cP^b$. 
\begin{Definition}
A $p$-sided {\it prism} in $\SLR$ space is an isometric image of a solid which is bounded by the side surfaces of a $p$-sided infinite prism
$\cP^i$ its base figur $\cP$ and the translated copy $\cP^t$ of $\cP$ by a fibre translation given by (1.5). 
\end{Definition}
The side faces $\cP$ and $\cP^t$ are called ,,{\it cover faces}"which are related by fibre translation along fibre lines joining their points.
\begin{Definition}
A $\SLR$ infinite prism is regular if $\cP^b$ is a regular $p$-gon with center at the origin in the ,,hyperbolic base plane" and 
the side surfaces are congruent to each other under an $\SLR$ isometry. 
\end{Definition}
\begin{Definition}
The regular $p$-sided {\it prism} in $\SLR$ space is a prism derived by the Definition 2.2 from a regular infinite prism (see Definition 2.3). 
\end{Definition}
\begin{Remark}
\begin{enumerate}
\item It is a natural assumption that the {\it ,,surfaces of the cover faces"} are derived as the images of the {\it ,,hyperbolic base plane"} at an isometry
of the $\SLR$ space i.e. the cover faces lie in Euclidean planes in the model.   
\item It is clear that there exist for all $p\in \mathbb{N}, \ (p\ge 3)$ $p$-gonal $\SLR$ prisms and also regular prisms (see Fig.~2, $\cP^b$ coincide with
$\cP$ and they are regular hyperbolic $p$-gons).   
\item All cross-sections of a prism ,,parallel" (the intersecting plane are generated by $\SLR$ fibre translations from the base plane) 
to the base faces are congruent. Prisms are named for their base, 
e.g. a prism with a pentagonal base is called a pentagonal prism (see Fig.~2). 
\end{enumerate}
\end{Remark}
A family of closed sets called tiles forms a tessellation or tiling of a space 
if their union is the whole space and every two distinct sets in the family have disjoint interiors.
A tiling is said to be monohedral if all of the tiles are congruent to each other.
At present the space is the $\SLR$ and the tiles are congruent {\it regular infinite or bounded prisms} (see Definition 2.2-3).
A tiling is called face-to-face if the intersection of any two tiles is either empty or a common face of both tiles otherwise it is 
non-face-to-face.

If the prisms are bounded then each vertex of a tiling is proper point of $\SLR$, 
thus the prism is a ,,$\SLR$ polyhedron" having at each vertex 
one ,,$p$-gonal cover face" (it is not absolutely polygon) and two skew ,,quadrangles" which lie on one-sheeted hyperboloid surfaces in the model. 
\begin{figure}[ht]
\centering
\includegraphics[width=6cm]{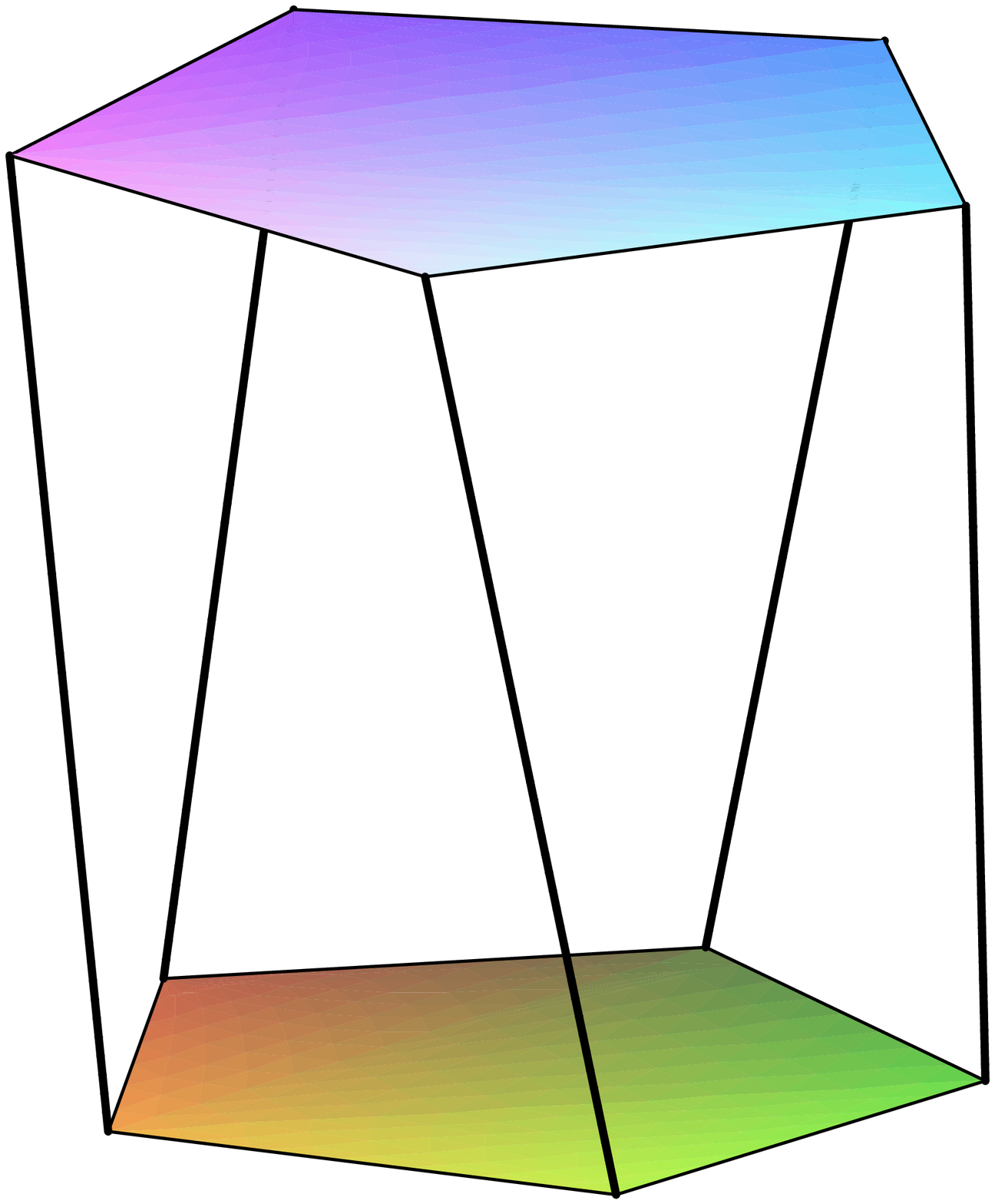} \includegraphics[width=6cm]{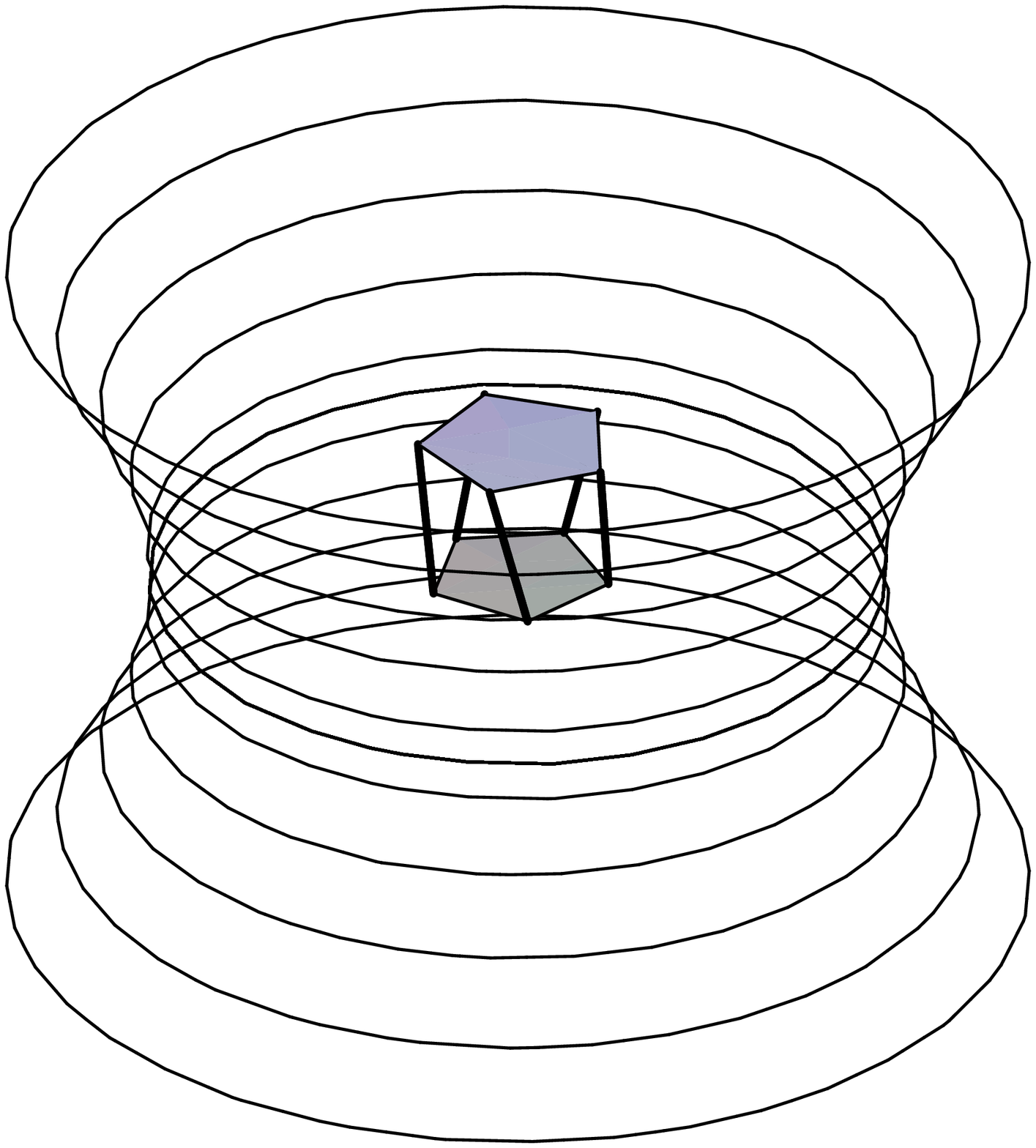}
\caption{Regular pentagonal prism}
\label{}
\end{figure}
\subsection{Regular infinite prism tilings}
First, we assume that $\cT_p^i(q)$ is a regular infinite prism tiling in the $\SLR$ space, which can be derived by a rotation subgroup 
$\bG_p^R(q)$ of the symmetry group $\bG_p(q)$ of $\cT_p^i(q)$.  $\bG_p^R(q)$ is generated by rotations $r_1;r_2;\dots;r_p$ with angle 
$\omega=\frac{2\pi}{q(p)}$ ($q \in \mathbb{N}$, $q(p)$ depends on the parameter $p$) about the fibre lines 
$f_1;f_2; \dots f_p$ through the vertices of the given $\SLR$ $p$-gon $\cP^b$ and let 
$\cP_p^i(q)$  be one of its tiles where we can suppose without loss of generality that its {\it $p$-gonal base figure} $\cP$ 
(and so $\cP^b$ as well) is centered at the origin.

The vertices $A_1A_2A_3 \dots A_p$ of the base figur coincide with the vertices of a regular
hyperbolic $p$-gon in the base plane with centre at the origin and we can introduce the following homogeneous coordinates to neigbouring vertices 
of the base figur of $\cP_p^i(q)$ in the hyperboloid model of $\widetilde{\cH}=\SLR$.
\begin{equation}
\begin{gathered}
A_1(1;0;0;x_3), \ A_2\Big(1;0;x_3\sin\Big(\frac{2\pi}{p}\Big);x_3\cos\Big(\frac{2\pi}{p}\Big)\Big), \\ 
A_3\Big(1;0;x_3\sin\Big(\frac{4\pi}{p}\Big);x_3\cos\Big(\frac{4\pi}{p}\Big)\Big). \tag{2.1}
\end{gathered}
\end{equation}
It is clear that the side curves $c(A_iA_{i+1})$ $(i=1\dots p, ~ A_{p+1} \equiv A_1)$ 
of the base figur are derived from each other by $\frac{2\pi}{p}$ rotation about the $x$ axis, so there are congruent in $\SLR$ sense.  
The necessary requirement to the existence of $\cT_p^i(q)$ that the surfaces of the neigbouring side faces of $\cP_p^i(q)$ are derived from each other 
by rotation with angle $\omega=\frac{2\pi}{q}$ \Big($ \frac{2p}{p-2} < q \in \mathbb{N}$\Big) 
about their ,,common fibre line".  

The isometry group of $\SLR$ leave invariant the hyperboloid $\cH$ and the fibre lines thus it is sufficient to consider the base $p$-gonal figur 
$A_1A_2A_3 \dots A_p$. Therefore, we have to require to the existence of a regular infinite $p$-gonal prism tiling $\cT_p^i(q)$ 
that the rotation $r_j(\omega)$
$(j=1,2,\dots,p)$ above the fibre lines $f_i$ (see (1.12)) maps the corresponding side face onto the neighbouring one: 
\begin{equation}
\begin{gathered}
r_1(\omega)~:~ [f_{p};f_1] \rightarrow [f_1;f_2], \ r_2(\omega)~:~ [f_1;f_2] \rightarrow [f_2;f_3], \\
r_3(\omega)~:~ [f_2;f_3] \rightarrow [f_3;f_4], \dots, r_p(\omega)~:~ [f_{p-1};f_p] \rightarrow [f_p;f_1]. \tag{2.2}
\end{gathered}
\end{equation}
\begin{Remark}
The isometries $r_i(\omega)$ $(i=1,2,\dots,p)$ map $\cP_p^i(q)$ onto its side face adjacent prisms, as well.   
\end{Remark}

$\cP_p^i(q)$ has rotational symmetry of the $2p$~th order about the $x$ axis therefore it is sufficient to require to the existence of $\cT_p^i(q)$
that e.g. $r_2(\omega)~:~ [f_1;f_2] \rightarrow [f_2;f_3]$. 
\begin{theorem}
There exist regular infinite prism tilings $\cT_p^i(q)$ for each $3 \le p \in \mathbb{N}$ where $q>\frac{2p}{p-2}$.    
\end{theorem}
{\bf Proof:} We have to prove two statements:
\begin{enumerate}
\item There are appropriate vertices (so ,,side fibre lines") of the base figur i.e. there is parameter $x_3$ so that $r_2(A_1)=A'_1$ 
lies on the fibre line through $A_3$. 
\item There are convenient side surfaces containing the corresponding side fibre lines i.e. there is a convenient side curve $c_{A_1A_2}$ 
of the base figur between $A_1$ and $A_2$ which image $c'_{A_1A_2}$ at rotation $r_2$ lies on the side surface generated by base side curve
$c_{A_2A_3}$.
\end{enumerate}
\begin{enumerate}
\item[({\bf i.})] We translate the points $A_1$, $A_2$, $A_3$ by $\SLR$ translation $\bT$ which map the point $A_2$ into the origin 
$$
\bT: A_1 \rightarrow A_1^T; \  \bT: A_2 \rightarrow O; \ \bT: A_3 \rightarrow A_3^T.
$$
The trace points of the fibres through $A_1^T$ and $A_3^T$ on the base plane are denoted by $A_1^{T_*}$ and $A_3^{T_*}$. 
To the existence of $\cT_p^i(q)$ the rotation about the fibre line $f_2$ with angle $\frac{2\pi}{q}$ 
has to map the fibre $f_1$ to $f_3$ thus the rotation about the $x$ axis with the above angle map the fibre $f_1^T$ to the fibre line $f_3^T$.
The $\SLR$ rotation about the $x$ axis in the hyperboloid model is the same as the Euclidean one therefore the points $A_1^{T_*}$ and $A_3^{T_*}$
lie in a circle in the hyperbolic base plane. Moreover, there is a $0< x_3 \in \mathbb{R}$ where the angle 
$A_1^{T_*} O A_3^{T_*}=\frac{2p}{q}$ ($q> \frac{2p}{p-2}$) because the angle of a hyperbolic $p$-gon is continuously changed in the intervall 
$(\frac{2p}{p-2},0)$ if $x_3\in (0,\infty)$. Therefore, the first statement is proved.  
\item[({\bf ii.})]
We have proved that there is $x_3$ that $r_2(A_1)=A'_1 \in f_3$. The trace point of $A'_1$ on the base plane is $A_3 \in f_3$.
Let $F\in f_3$ be the midpoint of the fibre segment $A'_1 A_3$ in $\SLR$ sense. The fibre lines through the points of $A_2F$ straight segment
form a side surface $S_{A_2A_3}$ (lying on a one-sheeted hyperboloid surface). $S_{A_2A_3}$ is a convenient side surface of $\cP_p^i(q)$ because
the curves $c_{A_1A_2}$ and $c'_{A_1A_2}$ are congruent therefore
the geodesic distances between the points $A_2,A_3$ and $A_2,A'_1$ are equal and so they are points of a geodesic ball centered at $A_2$, 
moreover the points $A_3$ and $A'_1$ lie in the fibre line $f_3$ and by the conditions of the fibre lines follows, 
that the further fibres (for example the fibre $f_0$ described in Fig.~3-4) through the points of the segment $A_2A'_1$ intersect the curves $c_{A_1A_2}$ and $c'_{A_1A_2}$, respectively (see Fig.3-4).  
\begin{figure}[ht]
\centering
\includegraphics[width=6cm]{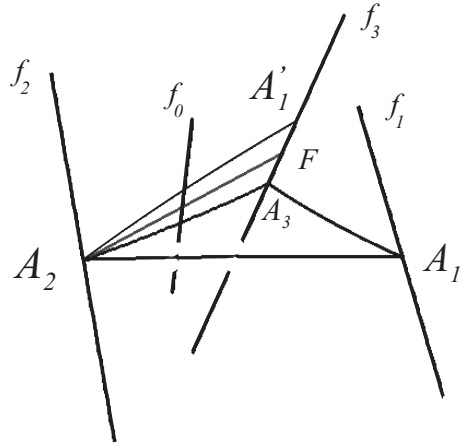}
\caption{The construction of $S_{A_2A_3}$ for regular infinite trigonal prism $\cP_3^i(7)$}
\label{}
\end{figure}
Therefore, the infinite (torus-like) prism tilings $\cT_p^i(q)$ exist. $\square$ 
\end{enumerate} 
\begin{Remark}
The equation of the curve $c_{A_1A_2}$ can be determined as the trace points (see (1.4) and (1.5)) of the fibres through the point of the segment $A_2F$.
The equations of the other side curves $c(A_iA_{i+1})$ $(i=2\dots p, \ A_{p+1} \equiv A_1)$ of the base figur are derived from the eqution of $c_{A_1A_2}$ 
by $\frac{2\pi}{p}$ rotation about $x$ axis (see Fig.~3 and Fig.~4).
\end{Remark}
\begin{figure}[ht]
\centering
\includegraphics[width=12cm]{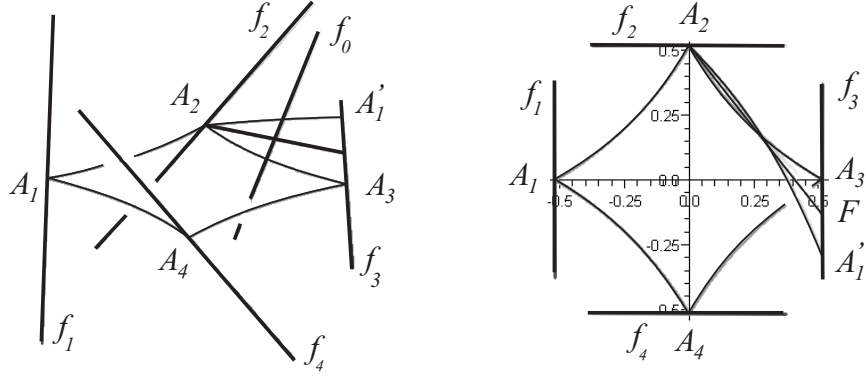}
\caption{The construction of $S_{A_2A_3}$ for regular infinite 4-gonal prism $\cP_4^i(6)$}
\label{}
\end{figure}
\subsubsection{Regular infinite trigonal prism tilings}
In this subsecton we determine the data of the existing (see Theorem 2.1) regular infinite trigonal prism tilings $\cT_3^i(q)$. 

The side faces of $\cP_3^i(q)$ are derived from each other 
by rotation with angle $\omega=\frac{2\pi}{q}$ ($ 6 < q \in \mathbb{N}$) 
about their ,,common fibre line". 

We use the homogeneous coordinates of vertices $A_1,A_2,A_3$ given in (2.1) depending on parameter $x_3$.
We have to determine parameter $x_3$  that the rotation $r_2(\omega)$
above fibre line $f_2$ (see (1.9)) maps the side face $[f_1;f_2]$ into the neighbouring one $[f_2;f_3]$.  

We obtain by above requirements an equation for the parameters $x_3$ and we get the following solution for each 
$7 \le q \in \mathbb{N}$: 

\begin{equation}
x_3=\sqrt{\frac{ \left( \sqrt {3}\cos \left( {\frac {2\pi }{q}} \right) -
\sin \left( {\frac {2\pi }{q}} \right)  \right) }{ \left( 2\,\sin
 \left( {\frac {2\pi }{q}} \right) +\sqrt {3} \right)}} \tag{2.3}
\end{equation}
Fig.~5 shows $\cP_3^i(7)$ with its base polygon. The equation of the curve $c_{A_1A_2}$ of $\cP_3^i(7)$ can be determined as the trace points (see (1.4) and (1.5)) 
of the fibres through the point of the segment $A_2F$ where $A'_3 \sim (1; 0.15072575; 0.23778592;$ $-0.18962794)$ and 
$F \sim (1; 0.07493964; 0.24918198; -0.16988939)$.
The equations of the other side curves $c(A_iA_{i+1})$ $(i=2,3, ~ A_4 \equiv A_1)$ of the base figur are derived from the equation of $c_{A_1A_2}$ 
by $\frac{2\pi}{3}$ rotation about $x$ axis (see Fig.~3 and Fig.~5). The data of $\cP_3^i(q)$ for some $\mathbb{N} \ni q > 6$ are collected in the Table 1.  
\begin{figure}[ht]
\centering
\includegraphics[width=6cm]{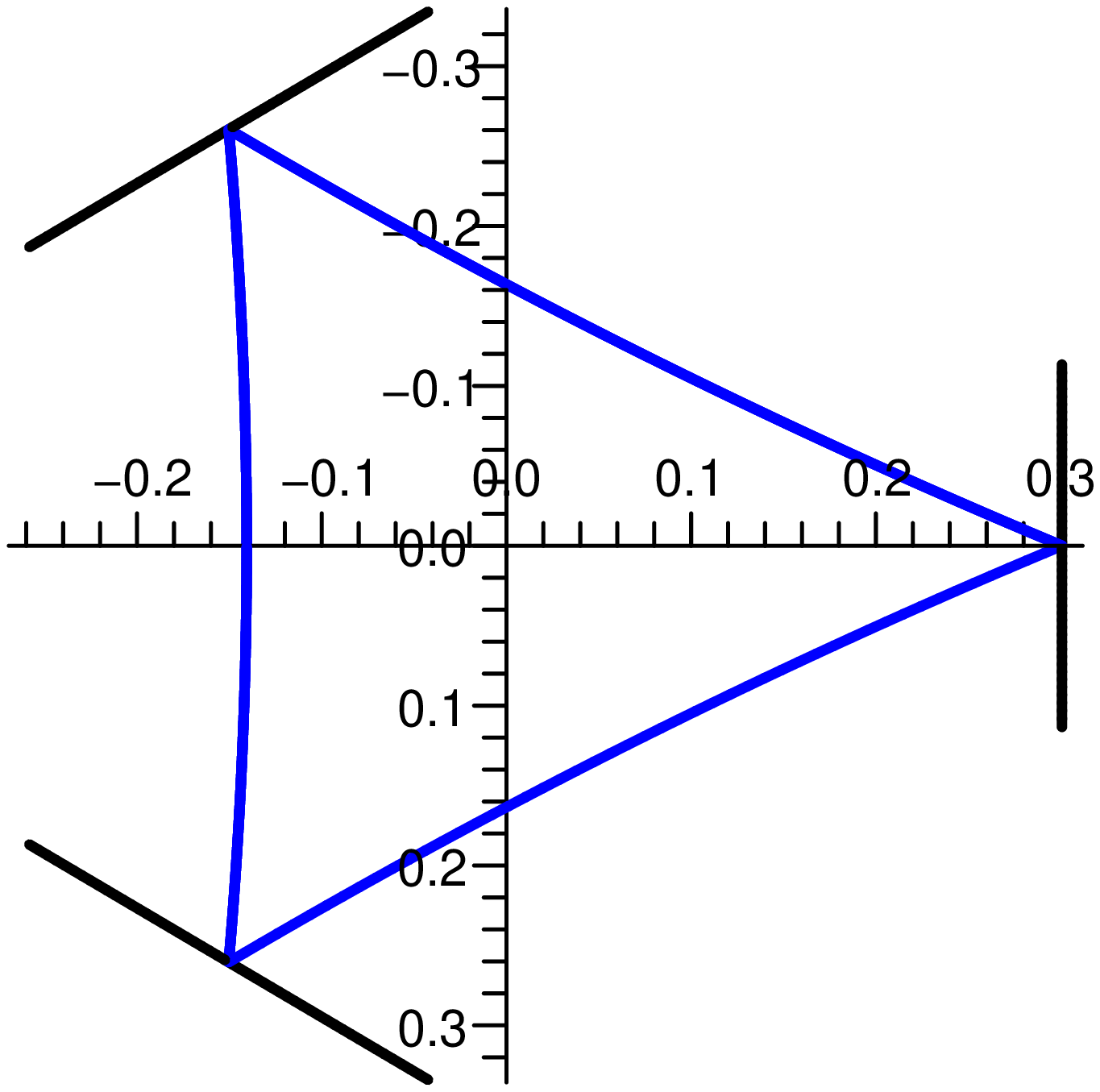} \includegraphics[width=7cm]{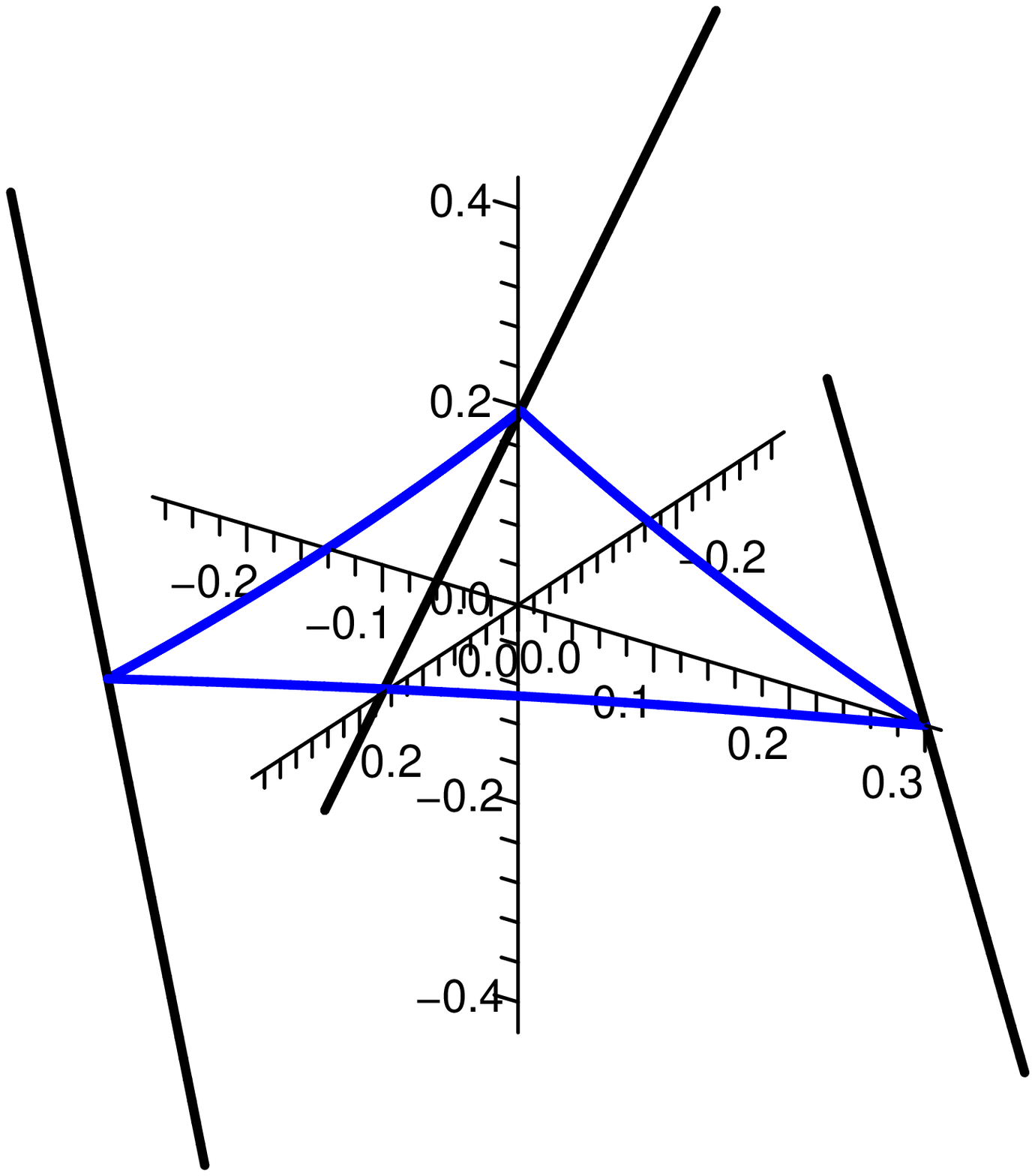}
\caption{Regular infinite trigonal prism $\cP_3^i(7)$ of $\cT_3^i(7)$}
\label{}
\end{figure}
\medbreak
\centerline{\vbox{
\halign{\strut\vrule\quad \hfil $#$ \hfil\quad\vrule
&\quad \hfil $#$ \hfil\quad\vrule
\cr
\noalign{\hrule}
\multispan2{\strut\vrule\hfill\bf Table 1 \hfill\vrule}%
\cr
\noalign{\hrule}
\noalign{\vskip2pt}
\noalign{\hrule}
(p, q)& x_3  \cr
\noalign{\hrule}
(3,7) & \approx 0.30007426  \cr
\noalign{\hrule}
(3,8) & \approx 0.40561640  \cr
\noalign{\hrule}
(3,9) & \approx 0.47611091  \cr
\noalign{\hrule}
(3,10) & \approx 0.50289355 \cr
\noalign{\hrule}
(3,50) & \approx 0.89636657  \cr
\noalign{\hrule}
(3,1000) & \approx 0.99457331  \cr
\noalign{\hrule}
}}}
\medbreak
We can determine the data of all regular infinite prism tilings $\cT_p^i(q)$ for given $3 \le p \in \mathbb{N}$ where $q>\frac{2p}{p-2}$.   
For example, we have described $\cP_4^i(6)$ with its base polygon in Fig~6, where the parameter $x_3=\frac{\sqrt{6}-\sqrt{2}}{2}$ . 
\begin{figure}[ht]
\centering
\includegraphics[width=6cm]{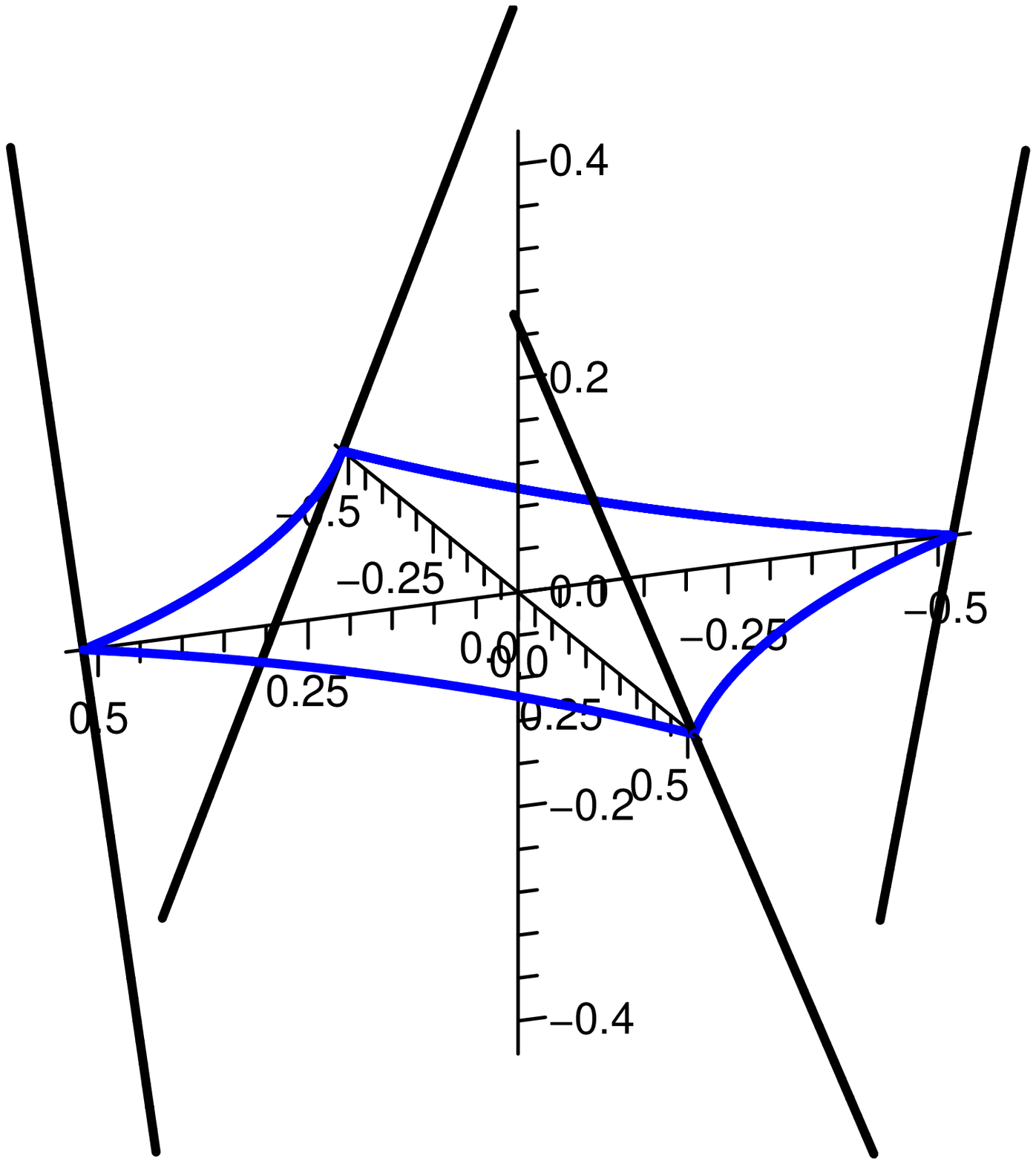} \includegraphics[width=7cm]{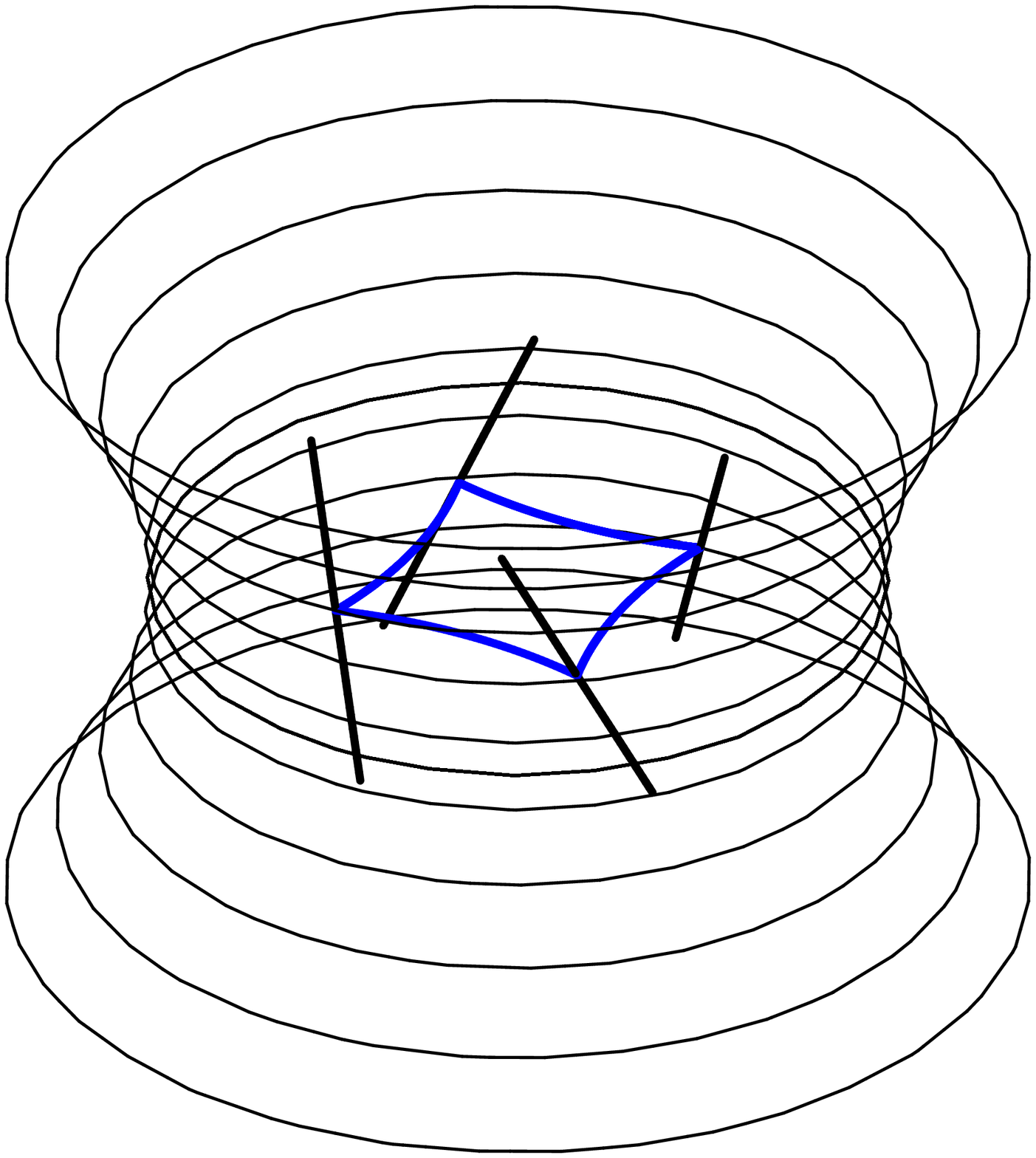}
\caption{Regular infinite 4-gonal prism $\cP_4^i(6)$ of infinite regular prism tiling $\cT_4^i(6)$}
\label{}
\end{figure}
\begin{figure}[ht]
\centering
\includegraphics[width=12cm]{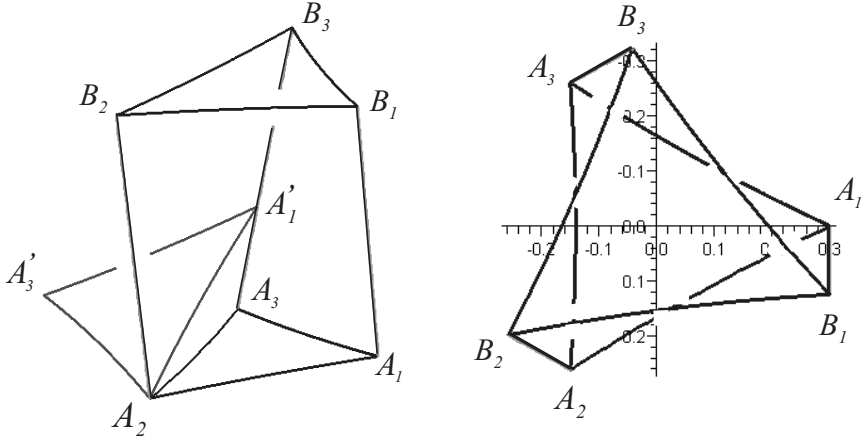}
\caption{Regular trigonal prism $\cP_3(7)~ (A_1A_2A_3B_1B_2B_3)$ with the base figur $A'_1A_2A'_3$ of its neighbouring prism.}
\label{}
\end{figure}
\subsection{Regular prism tilings}
In this section we study the regular (bounded) prism tilings in the $\SLR$ space. We can derive regular prism tilings $\cT_p(q)$ from the infinite regular
prism tilings $\cT_p^i(q)$ by the following way:
\begin{enumerate}
\item 
Let as suppose that $\cT_p^i(q)$ a regular infinite (or torus-like) prism tiling and let 
$\cP_p^i(q)$ be one of its tiles where $\cP$ (and so $\cP^b$ as well) is centered at the origin. 
Its {\it $p$-gonal base figure with vertices $A_1A_2A_3 \dots A_p$} in the hyperbolic base plane is derived as the intersection of 
$\cP_p^i(q)$ with the ,,base plane" of the model. It is clear that the side curves $c(A_iA_{i+1})$ $(i=1\dots p, ~ A_{p+1} \equiv A_1)$ 
of the base figur are derived from each other by $\frac{2\pi}{p}$ rotation about $x$ axis, so there are congruent in $\SLR$ sense.  
The corresponding vertices $B_1B_2B_3 \dots B_p$ are generated by a fibre translation $\tau$ given by (1.3) 
with parameter $\phi\in \mathbb{R}\setminus \{0\}$. The cover faces $A_1,\dots, A_p$, $B_1,\dots, B_p$ and the ,,side surfaces" form an 
$p$-sided regular prism $\cP_p(q)$ in $\SLR$.     
\item It is clear, that its images by the translation group $\langle \tau \rangle$ fill the regular infinite prism $\cP_p^i(q)$ without overlap.
\item $\cT_p^i(q)$ is generated by rotations $r_1;r_2;\dots;r_p$ with angle 
$\omega=\frac{2\pi}{q}$ \Big($ \frac{2p}{p-2} < q \in \mathbb{N}$\Big) about the fibre lines $f_1;f_2; \dots f_p$ through the vertices 
$A_1A_2A_3 \dots A_p$ therefore
we obtain a regular prism tiling $\cT_p(q)$, as well. 
\end{enumerate}
The images of the planes of equations $x=k$ ($k\in \mathbb{R}$) are invariant under rotations about the fibre line through the origin.
Therefore, their maps at an arbitrary translation, given by parameters $(t_0;t_1;t_2;t_3)$ (see (1.7)), are invariant planes 
under rotation $\bR_T(\omega)$ about the fibre line through the point $T(t_0;t_1;t_2;t_3)$ (see (1.9)). We get the next Lemma by (1.7). 
\begin{lemma}
The rotation $\bR_X(\omega)$ ($k\in \mathbb{R}$) leave invariant the planes of equations 
\begin{equation}
x(kt_1-t_0)+y(t_3-kt_2)-z(kt_3+t_2)+t_0k+t_1=0. \tag{2.4}
\end{equation}
\end{lemma}
Thus, the orbit of the point $A_1(1;0;0;x_3)$ lies by Lemma 2.2 at the rotation $r_2(\alpha)$ in the plane  
\begin{equation}
\begin{gathered}
S_2 \equiv -x+y\Big(x_3\cos\Big(\frac{2\pi}{p}\Big)-k x_3 \sin\Big(\frac{2\pi}{p}\Big)\Big)-\\ 
-z\Big(k x_3\cos\Big(\frac{2\pi}{p}\Big)+
x_3 \sin\Big(\frac{2\pi}{p}\Big)\Big)
+k=0, ~ \text{where} \ k=\frac{x_3^2 \sin\big(\frac{2\pi}{p}\big)}{1-x_3^2 \cos \big(\frac{2\pi}{p}\big)}. 
\end{gathered} \tag{2.5}
\end{equation}

It is clear, that the base plane and $S_2$ (see (2.5)) are different planes therefore 
the immediate consequence of the above Lemma 2.2 is the following
\begin{theorem}
There exist infinite many regular $p$-gonal non-face-to-face $\SLR$ prism tilings $\cT_p(q)$ for parameters $p \ge 3$ 
where $ \frac{2p}{p-2} < q \in \mathbb{N}$ but there is no face-to-face one.  
\end{theorem}
It is interesting to consider further tilings in the $3$-dimensional Thurston geometries, because important informations
of the ,,crystal structures" are included by the ,,space filling polyhedra". 

In this paper we have mentioned only some problems in discrete geometry of the $\SLR$ space, but we hope that from these 
it can be seen that our projective method suits to study and solve similar problems (see \cite{MSz12}, \cite{Sz07}, \cite{Sz10-1}, \cite{Sz12-1},
\cite{Sz12-2}). 



\begin{thebibliography}{MPSz98}
%
 \bibitem[1]{M97}
  {Moln{\'a}r,~E.}
  The projective interpretation of the eight 3-di\-men\-sional homogeneous geometries. 
  \emph{Beitr{\"a}ge zur Algebra und Geometrie (Contributions to Algebra and Geometry),}
 {\bf38} (1997) No.~2, 261--288.
  %
 \bibitem[2]{MSz}
 {Moln{\'a}r,~E.~--~Szirmai,~J.}
  Symmetries in the 8 homogeneous 3-geometries. 
  \emph{Symmetry: Culture and Science,}
  {\bf 21/1-3} (2010), 87-117.
 %
 \bibitem[3]{MSz12}
  {Moln{\'a}r,~E.~--~Szirmai,~J.}
   Classification of $\SOL$ lattices.
   \emph{Geometriae Dedicata,}
  (to appear) (2012), DOI: 10.1007/s10711-012-9705-5.
%
 \bibitem[4]{S}
 {Scott,~P.}
  The geometries of 3-manifolds.  
  \emph{Bull. London Math. Soc.}, {\bf15} (1983) 401--487. (Russian translation: Moscow "Mir" 1986.)
 %
 \bibitem[5]{Sz07}
 {Szirmai,~J.}
 The densest geodesic ball packing by a type of $\NIL$ lattices. 
 \emph{Beitr{\"a}ge zur Algebra und Geometrie (Contributions to Algebra and Geometry),}
 {\bf48(2)} (2007) 383--398.
 %
 \bibitem[6]{Sz10-1}
 {Szirmai,~J.}
 The densest translation ball packing by fundamental lattices in $\SOL$ space. 
 \emph{Beitr{\"a}ge zur Algebra und Geometrie (Contributions to Algebra and Geometry),}
 {\bf 51(2)} (2010), 353--373. 
 %
 \bibitem[7]{Sz10-2}
 {Szirmai,~J.}
 Geodesic ball packing in $\SXR$ space for generalized Coxeter space groups.
 \emph{Beitr{\"a}ge zur Algebra und Geometrie (Contributions to Algebra and Geometry),}
 {\bf 52(2)} (2011), 413--430. 
 %
 \bibitem[8]{Sz10-3}
 {Szirmai,~J.}
 Geodesic ball packing in $\HXR$ space for generalized Coxeter space groups.
 \emph{Mathematical Communications, to appear 2012}.
 %
 \bibitem[9]{Sz12-1}
 {Szirmai,~J.}
 Lattice-like translation ball packings in $\NIL$ space. 
 \emph{Publ. Math. Debrecen}
 {\bf 80/3-4} (2012), 427--440 (DOI: 10.5486/PMD.2012.5117). 
  %
 \bibitem[10]{Sz12-2}
 {Szirmai,~J.}
 On lattice coverings of the $\NIL$ space by congruent geodesic balls. 
 \emph{Mediterranean Journal of Mathematics}
 (to appear) [2012], DOI: 10.1007/s00009-012-0211-7. 
 %
 \bibitem[11]{T}
 {Thurston,~W.~P.} (and {\sc Levy,~S.} editor) 
 \emph{Three-Dimensional Geometry and Topology.}  
 {Princeton University Press,} Princeton, New Jersey, Vol {\bf 1} (1997).
 %
 \end{thebibliography}
\end{document}